\documentclass[12pt,draft]{iopart}

\usepackage{graphicx}
\usepackage{amssymb}
\usepackage{amsthm}
\usepackage{amsfonts}

\newtheorem{prop}{Proposition}[section]
\newtheorem{cor}[prop]{Corollary}

\newtheorem{thm}[prop]{Theorem}
\newtheorem{lemma}[prop]{Lemma}

\theoremstyle{definition}

\newtheorem{defn}[prop]{Definition}

\newcommand{\de}{\,\rmd}
\newcommand{\ee}{\rme}

\newcommand{\Rc}{\mathrm{Rc}}
\newcommand{\Div}{\mathrm{div}}
\providecommand{\text}[1]{{#1}}

\newcommand{\grad}{\nabla}
\newcommand{\p}{\partial}
\newcommand{\ga}{\alpha}
\newcommand{\N}{\nabla}
\newcommand{\abs}[1]{\ensuremath{\left|{#1}\right|}}
\newcommand{\norm}[1]{\ensuremath{\left\|{#1}\right\|}}
\newcommand{\ip}[2]{\left<{#1},{#2}\right>}

\newcommand{\gb}{\beta}
\newcommand{\gd}{\delta}
\newcommand{\gD}{\Delta}
\newcommand{\gf}{\varphi}
\newcommand{\gl}{\lambda}
\newcommand{\gS}{\Sigma}
\newcommand{\gO}{\Omega}
\newcommand{\gv}{\nu}
\newcommand{\gP}{\Pi}

\newcommand{\gm}{\mu}
\renewcommand{\ge}{\epsilon}

\renewcommand{\Re}{\ensuremath{\mathbb{R}}}

\newcommand{\ADM}{{\rm ADM}}
\newcommand{\Hk}{{\rm H}}
\newcommand{\ZAS}{{\rm ZAS}}
\newcommand{\Reg}{{\rm reg}}

\newcommand{\E}{\mathcal{E}}
\newcommand{\F}{\mathcal{F}}
\renewcommand{\bar}{\overline}
\renewcommand{\tilde}{\widetilde}
\newcommand{\bigo}{\mathcal{O}}

\begin{document}
\title{Zero area singularities in general relativity and inverse mean curvature
flow}
\author{Nicholas P Robbins}
\address{Gettysburg College, 300 North Washington Street, Gettysburg, PA 17325}
\ead{nrobbins@gettysburg.edu}

\begin{abstract}
First we restate the definition of a Zero Area
Singularity, recently introduced by H.L.\ Bray. We then consider several
definitions of mass for these singularities.  We use the Inverse Mean Curvature
Flow to prove some new results about the mass of a singularity, the ADM mass of
the manifold, and the capacity of the singularity. 
\end{abstract}
\ams{53, 83}
\submitto{\CQG}
\maketitle

\section{Introduction}

This paper consists of a study of some of the properties of Zero Area
Singularities, as recently introduced by Bray in \cite{Hugh-INI} and developed
by Bray and Jauregui in \cite{Hugh-Jeff}. The motivating example of which is the
spatial
Schwarzschild metric with a negative mass parameter: 
\begin{equation}
 g_{ij} = \left(1+\frac{m}{2r}\right)^4\gd_{ij}\qquad\qquad m<0.
\end{equation}
In addition to being historically and physically important, the
Schwarzschild solution is of particular mathematical interest since it
is the case of equality of the Riemannian Penrose conjecture
\cite{Hugh-Penrose}, and, in
the case when $m=0$, it is the case of equality of the Riemannian
Positive Mass Theorem \cite{SY-PMT}. Thus this metric, and its generalizations,
show promise as objects of study. For a further development of ZAS, as well as
an overview of some of the negative mass results in the field, see
\cite{Hugh-Jeff}.

The main results of this paper are, once we have defined the mass of a 
singularity, to extend the results of Huisken and Ilmanen in \cite{IMCF} 
to manifolds containing a single zero area singularity and a 
relationship between the capacity of the singularity and its mass.

\section{Definitions} 
\label{Defns}
\subsection{Asymptotically Flat Manifolds}
We will use the following definition of asymptotic flatness.
\begin{defn}(\cite{IMCF})
  A Riemannian 3-manifold $(M,g)$ is called \emph{asymptotically flat}
  if it is the union of a compact set $K$, and sets $E_i$, each
  diffeomorphic to the complement of a compact set $K_i$ in $\Re^3$,
  where the metric on each $E_i$ satisfies
\begin{equation}
  \abs{g_{ij}-\gd_{ij}}\leq \frac{C}{\abs{x}},\qquad
  \abs{g_{ij,k}}\leq \frac{C}{\abs{x}^2}
\end{equation}
as $\abs{x}\to\infty$. Derivatives are taken in the flat metric
$\gd_{ij}$ on $x\in\Re^3$. Furthermore the Ricci curvature must
satisfy 
\begin{equation}
  \Rc \geq -\frac{Cg}{\abs{x}^2}.
\end{equation}
The set $E_i$ is called an \emph{end} of $M$.
\end{defn}
A manifold may have several ends, but our results will be
relative to a single end. We will also be using the ADM mass of an asymptotically flat manifold and Hawking mass, capacity, and minimizing hull of a surface:
\begin{defn}\nocite{adm}
The \emph{ADM mass} of an end of an asymptotically flat manifold is 
\begin{equation}
  m_{\ADM}=\lim_{r\to\infty}\frac1{16\pi}\int_{S_\gd^r}\left(g_{ij,i}
    -g_{ii,j}\right)n^j \de\gm.
\end{equation}
\end{defn}

\begin{defn}
  The \emph{Hawking mass} of a surface $\gS$ is given by
  \begin{equation}
    m_{\Hk}=\sqrt{\frac{\abs{\gS}}{16\pi}}\left(1-\frac1{16\pi}\int_\gS
      H^2\right).
  \end{equation}
\end{defn}

\begin{defn}
  Let $\gS$ be surface in an asymptotically flat manifold
  $M$. Define the \emph{capacity} of $\gS$ by
  \begin{equation}
   C(\gS)= \inf\left\{\left. \int_M \norm{\nabla \gf}^2 \de V \right| 
        \gf(\gS)=1, \gf(\infty)=0\right\}.
  \end{equation}
\end{defn}
It is worth noting that if $\gS$ and $\gS'$ are two surfaces in $M$ so that
$\gS$ divides $M$ into two components, one containing infinity
and the other containing $\gS'$, then $C(\gS')\leq C(\gS)$
since the infimum is over a larger set of functions. 

\begin{defn}
Let $\gS$ be a surface that is the boundary of an open set, $E$, in a manifold $M$. 
We call $\gS$ a \emph{minimizing hull} if
\begin{equation}
    \abs{\p E \cap K} \leq \abs{\p F\cap K}
\end{equation}
for any $F$ containing $E$ where $K$ is a compact set containing $F\setminus E$.\end{defn}

\subsection{Definition and Mass of Zero Area Singularities}
The basic example of a Zero Area Singularity is the negative
Schwarzschild solution. This is the manifold $\Re^3\setminus B_{-m/2}$
with the metric
\begin{equation}
  g_{ij}  = \left(1+\frac{m}{2r}\right)^{4}\gd_{ij}
\end{equation}
where $m<0$. This manifold fails the requirements of the positive
mass theorem since it is not complete: geodesics reach the sphere at
$r=-m/2$ in finite distance. A straightforward calculation shows that
the ADM mass of this manifold is given by $m$.  Furthermore the far
field deflection of geodesics is the same as for a Newtonian mass of
$m$.  These results are identical to the same results for a positive
mass Schwarzschild solution. 

Two important aspects of this example will be incorporated into the
definition of a Zero Area Singularity. One is that the point
itself is not included. We
still must  describe the behavior of surfaces near the singularity. The
manifold in that region should have surfaces whose areas converge to
zero. In addition the capacity of these surfaces should go to zero.
The second aspect is the presence of a background metric, in this case
the flat metric.  This background metric will provide a location where
we can compute information about the singularity.

We have to be careful, since sometimes we are dealing with the toplogical
manifold $M$ and sometimes with the Riemannian manifold $M\setminus \gP$. So we
should clarify what we mean by ``convergence'' of surfaces. We again, follow
\cite{Hugh-Jeff} and restrict a surface to mean a $C^\infty$ closed embedded
2-manifold in the interior of of $M$ that is the boundary of a open region
$\gO$. We will mostly be concerned with surfaces converging to $\gP$. In this
case, for a surface sufficiently close to $\gP$, we can consider coordinates
$(x,s)$ on a tubular neighborhood of $\gP$, where $x\in \gP$ and $s\in[0,\ge)$.
Then we will restrict ourselves further to ``graphs'' over $\gP$. That is,
surfaces that can be written as $(x,s(x))$.  For such surfaces we define
convergence as follows.

\begin{defn}
Let $\{\gS_i\}$ be a sequence of surfaces that are graphs over $\gP$. Hence each
$\gS_i$ can be parametrized as $(x,s_i(x))$. We say that $\{\gS_i\}$
\emph{converges} to $\gP$ in $C^k$ if the functions $s_i:\gP\to(0,\ge)$ converge
to $0$ in $C^k$.   
\end{defn}
With this in hand we make the following definition.
\begin{defn}
  Let $M^3$ be a smooth manifold with boundary, where the boundary is
  compact.  Let $\gP$ be a compact connected component of the boundary
  of $M$. Let the interior of $M$ be a Riemannian manifold with smooth
  metric $g$. Suppose that, for any smooth family of surfaces, $\{\gS_i\}$,
  converging in $C^2$ to $\gP$,  
  the area of $\gS_i$ with respect to
  $g$ goes to zero as the surfaces converge to $\gP$.  Then $\gP$ is a
  \emph{Zero Area Singularity.}
\end{defn}
Bray and Jauregui in \cite{Hugh-Jeff} provide several equivalent conditions to
this definition.  We will use ZAS for the singular and plural of Zero Area
Singularity. A particularly useful class of these singularities are
\emph{regular}
Zero Area Singularities.
\begin{defn}
  Let $M^3$ be a smooth manifold with boundary. Let the boundary of
  $M$ consist of one compact component, $\gP$.  Let $\gP$ be a ZAS. If there is
  a smooth metric $\bar g$ on $M$
  and a smooth function $\bar \gf$ on $M$ with nonzero differential on
  $\gP$ so that $g=\bar\gf^4\bar g$, then we call $\gP$ a
  \emph{Regular} Zero Area Singularity.  We call the data $(
  M^3,\bar g, \bar \gf)$ a \emph{resolution} of $\gP$.
\end{defn}
Notice that while $\gP$ is topologically a surface, and it is a
surface in the Riemannian manifold $(M^3,\bar g)$, the areas of
surfaces near it in $(M^3\setminus\gP,g)$ approach zero, so we will
sometimes speak of $\gP$ as being a point $p$, when we are thinking in
terms of the metric $g$. Furthermore, notice that the requirement that
areas near $\gP$ go to zero under $g$ tells us that $\bar\gf=0$ on
$\gP$. For a regular ZAS, $g$ can be extended, as a symmetric two tensor, to
$\gP$.  In \cite{Hugh-Jeff} they consider local and global resolutions of ZAS.
However, since Geroch Monotonicity under IMCF requires our surfaces to be
connected, we only consider the case with a single ZAS. Thus we have no need for
a distinction between local and global resolutions.

We follow Bray in \cite{Hugh-INI} and define the mass of a regular ZAS as
follows:
\begin{defn}
  Let $(M^3,\bar g,\bar \gf)$ be a resolution of a regular ZAS $p=\gP$. Let
  $\bar\gv$ be the unit normal to
  $\gP$ in $\bar g$.  Then the
  \emph{regular mass} of $p$ is defined to be
  \begin{equation}
    m_{\Reg}(p)=-\frac14\left(\frac1\pi\int_{\gP}\bar\gv(\bar
      \gf)^{4/3}\, \bar{\de A}\right)^{3/2}.
  \end{equation}
\end{defn}
We also define the mass of a ZAS that may not be regular.
\begin{defn}\label{ZAS-mass-surf}
  Let $(M^3,g)$ be an asymptotically flat manifold, with a ZAS $p$.  
  Let $\gS_i$ be a smooth family of
  surfaces converging to $p$. Define
  $h_i$ by
  \begin{equation}
    \eqalign{\gD h_i = 0 \\
    \lim_{x\to\infty}h_i =1\\
    h_i = 0 \mbox{ on } \gS_i.}
  \end{equation}
  Then the manifold $(M,h^4_ig)$ has a ZAS
  at $\gS_i=p_i$ which is resolved by $(M,g,h_i)$.  Define the mass,
  $m_{\ZAS}(p)$, of $p$ to be
  \begin{equation} 
     \sup_{\{\gS_i\}}\limsup_{i\to \infty}
     -\frac14\left(\frac1\pi\int_{\gS_i} \gv(h_i)^{4/3}\, \de
      A\right)^{3/2}
      =\sup_{\{\gS_i\}}\limsup_{i\to \infty} m_{\Reg}(p_i).
  \end{equation}
  Here the outer sup is over all possible smooth families of surfaces
  $\{\gS_i\}$ which converge to $p$. 
\end{defn}

\section{Fundamental Results}
Before we continue we must verify that these definitions are
consistent.  First it must be verified that the regular mass of a
regular ZAS is indeed intrinsic to the singularity, as shown
in \cite{Hugh-INI}.
\begin{lemma}
  The regular mass of a  ZAS is
  independent of the resolution. 
  \begin{proof}
    Let $(M^3,\bar g, \bar \gf)$ and $(M^3,\tilde g, \tilde \gf)$ be
    two resolutions of the same ZAS, $p$. 
    Then define $\gl$ by
    $\bar \gf = \gl \tilde \gf$. Thus we note the following scalings:
    \begin{equation}
      \tilde g = \gl^4 \bar g\\
      \tilde{\de A}= \gl^4\bar{\de A}\\
      \tilde \gf = \gl^{-1}\bar\gf\\
      \tilde \gv = \gl^{-2}\bar\gv
    \end{equation}
    Now note that since $\tilde\gf,\bar\gf=0$ on $\tilde\gP,\bar\gP$, 
    \begin{equation}
      \tilde{\gv}(\tilde\gf)= \gl^{-2}\bar\gv\left(\gl^{-1}\bar\gf\right)=
      \gl^{-3}\bar\gv\left(\bar\gf\right)+\gl^{-4}\bar\gv(\gl)\bar\gf. 
    \end{equation}
    The last term, $\gl^{-4}\bar\gv(\gl)\bar\gf$, needs discussion.
    Both $\bar\gf$ and $\tilde\gf$ are smooth functions with zero set
    $\gP$ and they both have nonzero differential on $\gP$. Hence $\gl$
    is smooth. Thus since $\bar\gf$ goes to
    zero on $\gP$, this last term is zero on $\gP$.  Thus the mass of
    $p$ using the $(M^3,\tilde g,\tilde\gf)$ resolution is given by
    \begin{eqnarray}
      m_{\Reg}(p)=-\frac14\left(\frac1\pi\int_{\tilde\gP}\tilde\gv(\tilde
        \gf)^{4/3}\, \tilde{\de A}\right)^{3/2}\\
      =-\frac14\left(\frac1\pi\int_{\bar\gP}\left[\gl^{-2}
        \bar\gv(\gl^{-1}\bar\gf)\right]^{4/3}\, \gl^4\bar{\de A}\right)^{3/2}\\
      =-\frac14\left(\frac1\pi\int_{\bar\gP}\left[\gl^{-3}
        \bar\gv(\bar\gf)\right]^{4/3}\, \gl^4\bar{\de A}\right)^{3/2}\\
      =-\frac14\left(\frac1\pi\int_{\bar\gP}
        \bar\gv(\bar\gf)^{4/3}\,\bar{\de A}\right)^{3/2}.
  \end{eqnarray}
  \end{proof}
\end{lemma}
Definition~\ref{ZAS-mass-surf} seems to involve the entire manifold,
as the definition of $h_i$ takes place on the entire manifold. However
that isn't the case. The mass is actually local to the point $p$.
\begin{lemma}
  Let $(M^3,g)$ be a manifold with a ZAS
  $p$. Let $\tilde g$ be a second metric on $M$ that agrees with $g$
  in a neighborhood of $p$. Then the mass of $p$ in $(M^3,g)$ and
  $(M^3,\tilde g)$ are equal.
  \begin{proof}
    The goal is to show that for any selection of
    $\left\{\gS_i\right\}$, the series $ m_{\Reg}(p_i)$ and $\tilde
    m_{\Reg}(p_i)$ obtained in the calculation of the mass of $p$,
    with respect to $(M,g)$ and $(M,\tilde g)$ converge to the same value.  Let
    $S$ be a smooth, compact, connected surface separating $p$ from infinity
    and contained in the region where $g$ and $\tilde g$ agree.  Fix
    $i$ large enough so that $\gS_i$ is inside of $S$, and suppress
    the index $i$ on all our functions.  Then define the functions
    $h,\tilde h$ by
    \begin{equation}
      \eqalign{
      h=\tilde h = 0\mbox{ on } \gS_i\\
      \lim_{x\to\infty}h=\lim_{x\to\infty}\tilde h  = 1\\
      \gD h = \tilde \gD \tilde h  = 0.}
    \end{equation}
    Here $\gD$ and $\tilde \gD$ denote the Laplacian with respect to
    $g$ and $\tilde g$ respectively.

    Now inside $S$, $\gD=\tilde \gD$ since $g=\tilde g$. Thus there is
    only one notion of harmonic, and $h$ and $\tilde h$ differ only by
    their boundary values on $S$.  Let $\ge=1-\min_{S}\{h,\tilde h\}$.
    Consider the following two functions $f^-$ and $f^+$ defined between
    $S$ and $\gS_i$: \begin{equation}
      \eqalign{
      f^-=f^+=0\mbox{ on }\gS_i\\
      \gD f^- = \gD f^+ =0\\
      f^-=1-\ge\mbox{ on } S\\
      f^+=1\mbox{ on } S.}
    \end{equation}
    The maximum principle gives us, inside $S$,
    \begin{equation}
      f^+\geq h,\tilde h\geq f^-.
    \end{equation}
    Furthermore, since all four functions are zero on $\gS_i$,
    \begin{equation}
      \gv(f^+)\geq \gv(h),\gv(\tilde h)\geq \gv(f^-).
    \end{equation}
    Here $\gv$ is the normal derivative on $\gS_i$. Now define
    $\F(\gf)$ by the formula
    \begin{equation}
      \F(\gf) = \int_{\gS_i} \gv(\gf)^{4/3}\, dA.
    \end{equation}
    Then the ordering of the derivatives gives the ordering
     \begin{equation}
      \F(f^+)\geq \F(h),\F(\tilde h)\geq \F(f^-).\label{nrg-sand}
    \end{equation}
    However, since $f^-= (1-\ge)f^+$, 
    \begin{equation}
      \gv(f^-) = (1-\ge)\gv(f^+),
    \end{equation}
    hence,
    \begin{equation}
      \F(f^-) = (1-\ge)^{4/3} \F(f^+).
    \end{equation}
    Now, without loss of generality assume that the limit of the
    capacities of $\left\{\gS_i\right\}$ is zero, as the mass would be
    $-\infty$ otherwise.

  Thus as $i\to\infty$, $\gS_i$ has capacity going to zero (see
  subsection~\ref{cap-sub} for more discussion of capacity of ZAS.) Hence
$\ge_i$
  goes to zero, and so $\F(f^-_{i})/\F(f^+_i)$ goes to 1. Thus
  equation~\eref{nrg-sand} forces $\F(h_i)$ and $\F(\tilde h_i)$ to equality. 
  This forces the masses of $p_i$ in the two metrics to equality as well.
  \end{proof}
\end{lemma}
\begin{cor}
  In Definition~\ref{ZAS-mass-surf} we may replace the condition that
  $h_i$ be one at infinity with the condition that $h_i$ be one
  on a fixed surface outside $\gS_i$ for $i$ sufficiently large.
\end{cor}

\section{Zero Area Singularity Results}
\label{IMCF-ZAS}
\subsection{Zero Area Singularities and IMCF}
First recall that (weak) IMCF finds a (weak) solution to the equation:
\begin{equation}
  \Div_M\left(\frac{\grad u}{\abs{\grad u}}\right)=\abs{\grad u}.
\end{equation}
Whereever $u$ is smooth with $\grad u\neq 0$, the level sets of $u$ form a flow
of surfaces where the flow speed is given by $1/H$. 	

The main features we will be using of this are the following two facts:

\begin{thm}\label{asym-IMCF}
Assume $M$ is asymptotically flat, let $(N_t)_{t\geq t_0}$ be the surfaces
obtained from a weak solution to  IMCF in $M$. Then
\begin{equation}
  \lim_{t\to\infty} m_{\Hk}(N_t) \leq m_{\ADM}(M). 
\end{equation}
\end{thm}

\begin{thm}[Geroch Mononicity, 6.1 from \cite{IMCF}]
\label{G-mon}
Let $\tilde M$ be an asymptotically flat region of a manifold with $R>0$
exterior to a surface $\p \tilde M$.  For each connected component $N$  of
$\p\tilde M$, there exists a flow of compact $C^{1,\ga}$ surfaces $(N_t)_{t\geq
0}$, such that $N_0=N$, $m_\Hk(N_t)$ is monotone nondecreasing function for all
$t$ and for sufficiently large $t$, $N_t$ satisfies the IMCF.
\end{thm}
There is an assumption here that our starting surface is a minimizing hull. 
In order to apply IMCF to ZAS we will extend Geroch Monotonicity down to $t=0$
in the case where
our initial surface has negative Hawking mass.
\begin{lemma}\label{neg-mass-hull-lemma}
  Let $\gS$ be a surface in an asymptotically flat 3 manifold.  Let
  $\gS'$ be the boundary of the minimizing hull of $\gS$. Let $\gS$ or $\gS'$
  have negative Hawking mass. Then
  \begin{equation}
    m_\Hk(\gS)\leq m_\Hk(\gS').
  \end{equation}
\begin{proof}
  If $\gS'$ has nonnegative Hawking mass then $m_\Hk(\gS')\geq 0 \geq
  m_\Hk(\gS)$ and we are done. Thus we can assume that $\gS'$ has
  negative Hawking mass.  Since $\gS'$ has negative Hawking mass, it
  must intersect $\gS$ on a set of positive measure. Otherwise, $\gS'$
  would be a minimal surface, with Hawking mass
  $\sqrt{\frac{\abs{\gS'}}{16\pi}}>0$.  We define the following sets:
  \begin{equation}
    \gS_0 = \gS'\cap \gS\qquad
    \gS_+ = \gS'\setminus\gS_0\qquad
    \gS_- = \gS\setminus\gS_0
  \end{equation}
  Recalling that $\abs{\gS_+}\leq \abs{\gS_-}$ by the minimization
  property, and that $H=0$ on $\gS_+$, we observe the following:
  \begin{eqnarray}
    0>m_\Hk(\gS') &=
    \frac{\sqrt{\abs{\gS_0}+\abs{\gS_+}}}{(16\pi)^{3/2}}
    \left(16\pi-\int_{\gS_0}H^2\right)\\
    &\geq  \frac{\sqrt{\abs{\gS_0}+\abs{\gS_-}}}{(16\pi)^{3/2}}
    \left(16\pi-\int_{\gS_0}H^2\right)\\
    &\geq  \frac{\sqrt{\abs{\gS_0}+\abs{\gS_-}}}{(16\pi)^{3/2}}
    \left(16\pi-\int_{\gS_0}H^2-\int_{\gS_-}H^2\right)\\&=m_\Hk(\gS).\nonumber
  \end{eqnarray}
\end{proof}
\end{lemma}
With this lemma and Geroch Monotonicity we can prove  the following
lemma.
\begin{lemma}\label{geroch-2}
  Let $(M,g)$ be an asymptotically flat manifold with ADM mass $m$,
  nonnegative scalar curvature and a single regular ZAS $p$. Let $\{\gS_i\}$ be
  a smooth family of surfaces
  converging to $p$, which eventually have negative Hawking mass. Then
  for sufficiently large $i$,
   $ m_\Hk(\gS_i)\leq m$.
  \begin{proof}
    Since for large enough $i$, $\gS_i$ has non-positive Hawking mass, we
    can apply Lemma~\ref{neg-mass-hull-lemma} to show that $\gS_i'$
    must have larger Hawking mass. From this surface, we start Inverse
    Mean Curvature Flow. Geroch monotonicity tells us that the
    Hawking masses of the surfaces $N_t$ defined by IMCF starting with
    $\gS_i'$ only increase. Theorem 7.4 in \cite{IMCF} tells us that the
    increasing limit of the Hawking masses these surfaces is less than
    the ADM mass. Thus the Hawking mass of the starting surface was
    also less than the ADM mass.
    \end{proof}
\end{lemma}

Now we relate the limit of the Hawking masses to the regular mass.
\begin{lemma}\label{hawking-regular}
  Let $(M,g)$ be an asymptotically flat manifold with nonnegative
  scalar curvature and a single regular ZAS $p$.  Then there is a smooth family
  of surfaces $\{\gS_i\}$ converging to $p$ such that
\begin{equation}
  \lim_{i\to\infty} m_\Hk(\gS_i)=
  -\frac14\left(\frac1\pi\int_{\bar\gS} \bar\gv(\bar\gf)^{4/3}\bar{\de
      A}\right)^{3/2}=m_{\ZAS}(p).
\end{equation}
\begin{proof}
  The Hawking mass of a surface $\gS_i$ is given by
\begin{equation}
  m_\Hk(\gS_i) =
  \sqrt{\frac{\abs{\gS_i}}{16\pi}}\left(1-\frac1{16\pi}\int_{\gS_i}H^2\de
    A\right).
\end{equation}
Since the areas of the surfaces are converging to zero  we have
\begin{equation}
  \lim_{i\to\infty}  m_\Hk(\gS_i) =-  \lim_{i\to\infty}
  \frac{\sqrt{\abs{\gS_i}}}{(16\pi)^{3/2}}\int_{\gS_i}H^2\de A.
\end{equation}
By the H\"older inequality this is bounded as follows
\begin{equation}\label{loc-CS}
  - \frac{\sqrt{\abs{\gS_i}}}{(16\pi)^{3/2}}\int_{\gS_i}H^2\de A\leq -
  \frac1{(16\pi)^{3/2}}\left(\int_{\gS_i} H^{4/3}
    \de A\right)^{3/2}.
\end{equation}
Switching to the resolution space, we use the formula
\begin{equation}
  H =\bar\gf^{-2} \bar H + 4\bar\gf^{-3}\bar\gv(\bar\gf ).
\end{equation}
Putting this into the previous equation we get
\begin{eqnarray}
   \int_{\gS_i} H^{4/3}
    \de A&=\int_{\bar{\gS}_i} \left( \bar\gf^{-2} \bar H +4
      \bar\gf^{-3}\bar\gv(\bar\gf ) \right)^{4/3}\gf^4\bar{\de A}\\
    &=\int_{\bar{\gS}_i} \left( \bar\gf \bar H +4
      \bar\gv(\bar\gf ) \right)^{4/3}\bar{\de A}.
\end{eqnarray}
Since $\bar\gf$ is zero on $\bar\gS$ and $\bar H$ is bounded, the first
term goes to zero. The second term converges since the family of
surfaces $\{\gS_i\}$ are converging smoothly.
\begin{equation}
  \lim_{i\to\infty}\int_{\bar{\gS}_i} \left( \bar\gf \bar H +4
    \bar\gv(\bar\gf ) \right)^{4/3}\bar{\de A} =4^{4/3}\int_{\bar\gS}
  \bar\gv(\bar\gf )^{4/3}\bar{\de A}.
\end{equation}
Combining all of these equations we have
\begin{equation}
  \lim_{i\to\infty} m_\Hk(\gS_i)\leq
  -\frac14\left(\frac1\pi\int_{\bar\gS} \bar\gv(\bar\gf)^{4/3}\bar{\de
      A}\right)^{3/2}=m_{\ZAS}(p).\label{loc-HP}
\end{equation}
To see when this estimate is sharp, we look at inequality
\eref{loc-CS} since that is the only inequality is our estimate. In
the limit, this inequality is an equality exactly when the ratio of
the maximum and minimum values of $H$ approaches $1$. We choose a
resolution such that $\bar\gv(\bar\gf)=1$ on the boundary. We also
choose a family of surfaces $\gS_i$ given by level sets of $\bar\gf$.
Then if we look at the ratio
\begin{equation}
  \lim_{\gf\to0} \frac{H_{\mathrm{min}}}{H_{\mathrm{max}}}=
  \lim_{\gf\to0}\frac{\bar\gf \bar
    H_{\mathrm{min}}+4\bar\gv(\bar\gf)}{\bar\gf \bar
    H_{\mathrm{max}}+4\bar\gv(\bar\gf)},
\end{equation}
and remember that $\bar H$ is bounded, we see that the
$\bar\gv(\bar\gf)$ terms dominate, and as $\bar\gf\to0$, this ratio
approaches 1. Thus with this resolution and this family of surfaces,
inequality \eref{loc-HP} will turn to an equality.
\end{proof}
\end{lemma}
With these results we can prove
the following theorem
\begin{thm}\label{penrose-one-regular}
  Let $(M,g)$ be an asymptotically flat manifold with nonnegative
  scalar curvature and a single regular ZAS $p$.  Then the ADM mass of $M$ is at
least the mass of
  $p$.
\begin{proof}
  First consider the case when $p$ can be enclosed by a surface,
  $\gS$, with nonnegative Hawking mass. The minimizing hull of a
  surface with nonnegative Hawking mass has nonnegative Hawking mass.
  Thus we can run IMCF from $\gS'$, and the AMD mass of $M$ is at
  least $m_\Hk(\gS')\geq0$.  However, the regular mass of $p$ is
  always nonpositive so in this case we are done.

  Now assume that $p$ cannot be enclosed by a surface with nonnegative
  Hawking mass. By Lemma~\ref{geroch-2} we know that the ADM mass is
  greater than the Hawking masses of any sequence of surface
  converging to $p$ which have negative Hawking mass. By
  Lemma~\ref{hawking-regular} we know that there is a family of
  surfaces converging to $p$ which have the mass of $p$ as the limit
  of their Hawking mass, hence the ADM mass is greater then their
  Hawking masses which limit to the regular mass.
\end{proof}
\end{thm}

This can be extended to a general ZAS. However, first we need to consider the
effect of
multiplication by a harmonic conformal factor on the ADM mass of a
manifold. 
\begin{lemma}\label{harm-mass-mod}
  Let $(M^3,g)$ be an asymptotically flat manifold. Let $\gf$ be a
  harmonic function with respect to $g$ with asymptotic expansion
  \begin{equation}
    \gf = 1+\frac{C}{\abs{x}_g}+\bigo\left(\frac{1}{\abs{x}^2_g}\right).
  \end{equation}
  Then, if the ADM mass of
  $(M^3,g)$ is $m$, the ADM mass of $(M^3,\gf^4g)$ is $m+2C$. 
  \begin{proof}
    This is a direct calculation. We write $g^\gf = \gf^4g$, and
    calculate, only keeping the terms of lowest order in $\abs{x}^{-1}$
    since we are taking limits as $\abs{x}\to\infty$.
\begin{eqnarray}
  m_\gf&=
  \lim_{\abs{x}\to\infty}\frac1{16\pi}\int_{S^\gd}\left(g_{ij,i}^\gf
    -g_{ii,j}^\gf\right)n^j\de A\\
  &=
  \lim_{\abs{x}\to\infty}\frac{\gf^4}{16\pi}\int_{S^\gd}\left(g_{ij,i}
    -g_{ii,j}\right)n^j\de A\nonumber\\
&\qquad\qquad+
  \lim_{\abs{x}\to\infty}\frac{\gf^3}{4\pi}\int_{S^\gd}\left(\gd_{ij}\gf_i
    -\gd_{ii}\gf_j\right)n^j\de A\\
  &= \lim_{\abs{x}\to\infty}\gf^4m -
  \lim_{\abs{x}\to\infty}\gf^3\lim_{\abs{x}\to\infty}
  \frac{1}{4\pi}\int_{S^\gd}\left(\gf_j-3\gf_j\right)n^j\de A\\
  &=m -
  \lim_{\abs{x}\to\infty}\frac{1}{2\pi}\int_{S^\gd}\gf_jn^j\de A\\
  &=m - \lim_{\abs{x}\to\infty}\frac{1}{2\pi}\int_{S^\gd}\ip{\grad\gf}{\gv}\de
A\\
  &=m+2C.
\end{eqnarray}

\end{proof}
\end{lemma}

Using this we can now extend Theorem \ref{penrose-one-regular} to a
general ZAS. 
\begin{thm}\label{penrose-one-general}
  Let $(M,g)$ be an asymptotically flat manifold with nonnegative
  scalar curvature and a single zero area singularity $p$.
  Then $m$, the ADM mass of $M$, is at least the mass of $p$.
  \begin{proof}
    If the capacity of $p$ is nonzero, then the statement is trivial.
    Thus we assume the capacity of $p$ is zero.  Using the terminology
    of Definition \ref{ZAS-mass-surf}, Theorem
    \ref{penrose-one-regular} tells us that the ADM mass of
    $(M,h^4_ig)$ is at least the mass of the regular singularity at
    $\gS_i=p_i$. Each $h_i$ is defined by the equations
  \begin{equation}
     \eqalign{     \gD h_i = 0 \\
      \lim_{x\to\infty}h_i =1\\
      h_i = 0 \mbox{ on } \gS_i .}
    \end{equation}
    Thus it has asymptotic expansion
    \begin{equation}
      h_i = 1-\frac{C_i}{\abs{x}}+\bigo\left(\frac1{\abs{x}^2}\right).
    \end{equation}
    Where $4\pi C_i$ is the capacity of $\gS_i$.  Thus, the ADM mass,
    $m_i$, of $(M,h^4_ig)$ is given by $m-2C_i$. Now we know that
    $m_i\geq m_{\Reg}(p_i)$. Taking $\limsup$ of both sides gives
    us
    \begin{equation}
      \limsup_{i\to\infty}m_i\geq \limsup_{i\to\infty} m_{\Reg}(p_i)
    \end{equation}
    Since $C_i$ is going to zero, the left hand side is simply $m$,
    and so has no dependence on which $\{\gS_i\}$ we chose in our mass
    calculation. Thus we get
    \begin{equation}
      m \geq \sup_{\{\gS_i\}}\limsup_{i\to \infty} m_{\Reg}(p_i)=m_{\ZAS}(p).
    \end{equation}
    as desired.
  \end{proof}
\end{thm}

\subsection{Capacity and ZAS}\label{cap-sub}
The capacity of a surface provides a measure of its size as seen from
infinity. We extend the definition of the capacity of surface to the
capacity of a zero area singularity.  We then show that if a
ZAS has non-zero capacity the Hawking mass of any family of surfaces
converging to it must go to negative infinity. We now define the capacity of a singular point. The natural
definition is the one we want.
\begin{defn}
  Let $p$ be singular point in an asymptotically flat manifold $M$.
  Chose a sequence of surfaces $\gS_i$ of
  decreasing diameter enclosing $p$. Then define the \emph{capacity}
  of $p$ by the limit of the capacities of $\gS_i$.
\end{defn}
Before using this definition we have to show that it is well defined.
\begin{lemma}
  Let $\gS_i$ and $\widetilde\gS_i$ be two sequences of surfaces
  approaching the point $p$. If $\lim C(\gS_i)=K$, $\lim
  C(\widetilde\gS_i)=K$. Hence $C(p)$ is well defined.
\begin{proof}
  Since the $\gS_i$ are going to $p$, for any given $\widetilde
  \gS_{\tilde{i_0}}$, we can choose $i_0$ such that for all $i> i_0$, $\gS_i$
  is contained within $\widetilde \gS_{\tilde{i_0}}$. Thus if $\gf$ is a
  capacity test function for $\widetilde \gS_{\tilde{i_0}}$, i.e.\
  $\gf(\widetilde \gS_{\tilde{i_0}})=1$ and $\gf\to0$ at infinity, then $\gf$ is
  also a capacity test function for $\gS_i$. Since $C(\gS_i)$ is taking the
  infimum over a larger set of test functions than $C(\widetilde
  \gS_{\tilde{i_0}})$, $C(\gS_i)\leq C(\widetilde \gS_{\tilde{i_0}})$. Thus if
  we create a new sequence of surfaces $\bar \gS_i$, alternately choosing
  from $\gS_i$ and $\widetilde\gS_i$, such that each surface contains
  the next we get a nonincreasing sequence of capacities. Thus if
  either original sequence of surfaces has a limit of capacity, then
  this new sequence must as well, and it must be the same. Hence,
  $\lim_{i\to\infty}C(\gS_i)=\lim_{i\to\infty}C(\widetilde\gS_i)$.
\end{proof}
\end{lemma}

Now we look at the relationship between capacity and the Hawking mass
of a surface. We will use techniques similar to those used in
\cite{Bray-Miao}.

\begin{thm}
  \label{cap-mass-lemma}
  Let $M$ be an asymptotically flat 3 manifold with nonnegative scalar
  curvature, and ZAS $p$.  Let $\gS_i$ be
  a family of surfaces converging in $C^2$ to $p$.  Assume each $\gS_i$ is a
  minimizing hull. Assume the areas of $\gS_i$ are going to zero. Then
  if the Hawking mass of the surfaces is bounded below, the capacities
  of surfaces converging to $p$ must go to zero.
  \begin{proof}
    To use Geroch monotonicity, we need to know that our IMCF surfaces
    stay connected. In the weak formulation of IMCF, the level sets
    $\gS_t$ always bound a region in $\bar M$.  Thus if $\gS_t$ is not
    connected, one of its components $\gS_t^*$ must not bound a
    region. That is, $\gS_t^*$ is not homotopic to a point in $M$.
    Since $M$ is smooth, it must have finite topology on any bounded
    set. Thus we know that near $p$, there is a minimum size for a
    surface that does not bound a region. Call this size
    $A_\mathrm{min}$.  Thus if we have any surface that does not bound a
    region, it must have area greater then $A_{\mathrm{min}}$. The area of our
    surfaces grow
    exponentially. Thus if we restrict ourselves to starting IMCF with
    a surface with area $A_{\rm min}/\ee$, and only run the flow for
    time $1$, we will stay connected. At first glance it seems we may
    need to worry about the jumps in weak IMCF, however Geroch
    monotonicity doesn't depend on smoothness of the flow, and neither
    does the area growth formula. Thus even with jumps, the area of our surfaces
    will remain below $A_{\rm min}$.

    Now recall that capacity of a surface is defined by
    \begin{equation}
      C(\gS)= \inf\left\{\left. \int_M \norm{\nabla \gf}^2 \de V \right| 
        \gf(\gS)=1, \gf(\infty)=0\right\}.
    \end{equation}
    Here, the integral is only over the portion of $M$ outside of
    $\gS$.  Call this integral, $\E(\gf)$, the \emph{energy} of $\gf$.
    Thus for any $\gf$ with $\gf(\infty)=0$ and $\gf(\gS)=1$ we have
    $\E(\gf)\geq C(\gS)$. So we will find an estimate that relates the
    Hawking mass and the energy of a test function $\gf$.

    Choose a starting surface $\gS$ with sufficiently small starting
    area. Let $f$ be the level set function of the associated weak
    IMCF starting with the surface $\gS$. Call the resulting level
    sets $\gS_t$. Now if we use a test function of the form
    $\gf=u(f)$, then the energy of $\gf$ is given by
    \begin{equation}
      \E(\gf)= \int_M \norm{\N f}^2(u')^2 \de V.
    \end{equation}
    Since $f$ is given by IMCF, we know that $\norm{\N f} = H$ where
    $H$ is the mean curvature of the level sets. Next we use the
    co-area formula with the foliation $\gS_t$ and our integral becomes
    \begin{equation}
      \E(\gf)= \int_0^\infty (u'(t))^2\int_{\gS_t} \abs{H} \de A_t\,dt.
    \end{equation}
    Here the co-area gradient term cancels one of the $\abs{H}=\norm{\N f}$
    terms. Now we will bound the interior integral of curvature. We know that
    IMCF causes the Hawking mass to be nondecreasing in $t$. We first
    rewrite the definition of the Hawking mass $m_\Hk(\gS_t^i) = m(t)$
    as:
    \begin{equation}
      \int H^2 \de A_t = 16\pi\left(1 - m(t)\sqrt{\frac{16\pi}{A(t)}}\right).
    \end{equation}
    Here $A(t)$ is the area of $\gS_t$.  Since the Hawking mass is
    nondecreasing under IMCF we have:
    \begin{equation}
      \int H^2 \de A_t \leq 16\pi\left(1 - m(0)\sqrt{\frac{16\pi}{A(t)}}\right).
    \end{equation}
    Thus we can use Cauchy-Schwartz to get:
    \begin{equation}
      \int \abs H \de A_t \leq \sqrt{A(t)}\sqrt{ 16\pi\left(1 -
          m(0)\sqrt{\frac{16\pi}{A(t) }}\right)}.
    \end{equation}
    We can rewrite this as:
    \begin{equation}
      \int \abs H \de A_t \leq \sqrt{\ga A(t)+\gb\sqrt{A(t)}}.
    \end{equation}
    Furthermore, since $A(t)$ grows exponentially in $t$, we can write
    this as:
    \begin{equation}
      \int \abs H \de A_t \leq \sqrt{\ga \ee^t+\gb \ee^{t/2}}=v(t).
    \end{equation}
    Where $A_0$ has been absorbed into $\ga$ and $\gb$. Thus our
    energy formula has become
    \begin{equation}
      \E(\gf) \leq \int_0^\infty (u'(t))^2v(t)\de t.
    \end{equation}
    with
    \begin{equation}
      v(t)=\sqrt{\ga \ee^t+\gb \ee^{t/2}}
    \end{equation}
    where $\ga= 16\pi A_0$, $\gb = (16\pi)^{3/2}A_0^{1/2}\abs{m_0}$,
    and $A_0$ is $A(\gS_0)$. This means we can pick our test
    function $u(t)$ to be as simple as:
    \begin{equation}
      u(t)=\cases{
      1-t & $0\leq t\leq 1$ \\
      0&$t\geq 1$.\\}
\end{equation}
    Then our integral becomes:
    \begin{eqnarray}
      E(\gf)&\leq \int_0^1 v(t)\de t\\
    &\leq 2\sqrt\ga+2\sqrt\gb.
    \end{eqnarray}
    Since $m_{\Hk}(\gS)\leq \sqrt{\frac{\abs{\gS}}{16\pi}}$, $m_0$ is
    bounded above. By assumption $m_0$ is bounded below, so $\ga$ and
    $\gb$ are bounded by multiples of $A_0$ and $\sqrt{A_0}$
    respectively. 
  Thus $\E(\gf)$ goes to zero if $A_0\to 0$ and $m_0$
  is bounded.  Hence $C(p)$ must be zero since it is the infimum over
  a positive set with elements approaching zero.
  \end{proof}
  
\end{thm}

\begin{thm}[Capacity Theorem]
  \label{cap-mass-theorem}
  Let $M$ be an asymptotically flat 3 manifold with nonnegative scalar
  curvature, and ZAS $p$, such that there
  exists a family of surfaces, $\gS_i$, converging in $C^2$ to $p$.  
  Then if the capacity of $p$ is nonzero, the Hawking
  masses of the surfaces $\gS_i$ must go to $-\infty$.
  \begin{proof}
    Any such family of surfaces will generate a family, $\{\gS_i'\}$, of
    minimizing hulls that will also converge to $p$.  By Theorem
    \ref{cap-mass-lemma}, the masses of $\{\gS_i'\}$ must go to
    $-\infty$. Thus the masses of $\{\gS_i'\}$ must go to $-\infty$.
    Thus for sufficiently large $i$, the masses of the minimizing
    hulls are all negative.  From then on Lemma
    \ref{neg-mass-hull-lemma} applies, and the masses of $\gS_i$ must
    be less then the masses of $\gS_i'$. Hence they also converge to
    $-\infty$.
  \end{proof}

\end{thm}

\ack
The author of this paper was supported in part by NSF grant \#DMS-05-33551,
through Duke University.
\section*{References}
\bibliographystyle{unsrt}
\bibliography{ZAS_robbins}

\begin{thebibliography}{1}

\bibitem{Hugh-INI}
Hubert~L. Bray.
\newblock Negative point mass singularities in general relativity.
\newblock In {\em Global problems in mathematical general relativity}. Issac
  Newton Institute, University of Cambridge, August 2005.
\newblock Available at
  \texttt{http://www.newton.cam.ac.uk/webseminars/pg+ws/2005/gmr/0830/bray/}.

\bibitem{Hugh-Jeff}
Hubert~L. Bray and Jeffrey~L. Jauregui.
\newblock A geometric theory of zero area singularities in general relativity,
  2009.
\newblock Available at \texttt{http://arxiv.org/abs/0909.0522}.

\bibitem{Hugh-Penrose}
Hubert~L. Bray.
\newblock Proof of the {R}iemannian {P}enrose inequality using the positive
  mass theorem.
\newblock {\em Journal of Differential Geometry}, 59:177--2677, 2001.

\bibitem{SY-PMT}
R.~Schoen and S.-T. Yau.
\newblock Proof of the positive mass theorem {II}.
\newblock {\em Communications in Mathematical Physics}, 79:231--260, 1981.

\bibitem{IMCF}
Gerhard Huisken and Tom Ilmanen.
\newblock The inverse mean curvature flow and the {R}iemannian {P}enrose
  inequality.
\newblock {\em J. Differential Geom.}, 59(3):353--437, 2001.

\bibitem{adm}
R.~Arnowitt, S.~Deser, and C.~W. Misner.
\newblock Coordinate invariance and energy expressions in general relativity.
\newblock {\em Phys. Rev. (2)}, 122:997--1006, 1961.

\bibitem{Bray-Miao}
Hubert Bray and Pengzi Miao.
\newblock On the capacity of surfaces in manifolds with nonnegative scalar
  curvature.
\newblock {\em Invent. Math.}, 172(3):459--475, 2008.

\end{thebibliography}
\end{document}